\documentclass[12pt,a4paper]{article}
\usepackage[T1]{fontenc}
\usepackage[utf8]{inputenc}
\usepackage[british]{babel}
\usepackage{amsmath,amssymb,amsthm}
\usepackage{tensor}
\usepackage{xcolor}
\usepackage[all]{xy}
\usepackage{cite}
\usepackage{hyperref}
\usepackage{tikz}
\usetikzlibrary{knots}

\usepackage[shortlabels]{enumitem}
\setlist[enumerate]{leftmargin=*}
\setlist[itemize]{leftmargin=*}

\newcommand{\Z}{\mathbb{Z}}

\DeclareMathOperator{\Soc}{Soc}
\DeclareMathOperator{\Ker}{Ker}
\DeclareMathOperator{\Aut}{Aut}

\DeclareMathOperator{\Fix}{Fix}

\newtheorem{teo}{Theorem}[section]
\newtheorem{prop}[teo]{Proposition}

\newtheorem{lemma}[teo]{Lemma}
\newtheorem{cor}[teo]{Corollary}

\newtheorem{athm}{Theorem}

\theoremstyle{definition}

\theoremstyle{remark}

\theoremstyle{definition}
\newtheorem{ex}[teo]{Example}

\title{On complemented skew braces}
\date{}
\author{Yangming Li\thanks{School of Mathematics, Guangdong University of Education, Guangzhou, 510310, China \texttt{3355800852@qq.com}}
\and
 P. Pérez-Altarriba\thanks{Departament de Matem\`atiques, Universitat de Val\`encia, Dr.\ Moliner, 50, 46100 Burjassot, Val\`encia, Spain; \texttt{Pedro.A.Perez@uv.es}, \texttt{Adolfo.Ballester@uv.es}}
\and 
L. A. Kurdachenko\thanks{Department of Algebra and Geometry, Oles Honchar Dnipro National University, Dnipro 49010, Ukraine; \texttt{lkurdachenko@gmail.com}}
\and
A. Ballester-Bolinches\addtocounter{footnote}{-2}\footnotemark
}

\begin{document}
\maketitle

\begin{abstract}
We investigate complemented skew left braces, that is, skew braces in which every skew subbrace admits a complement. We establish that if every minimal subbrace of a finite skew brace is complemented, then the brace is supersoluble and all its skew subbraces of prime power order are elementary abelian. We also show that the converse does not hold in general, providing a finite counterexample. Furthermore, we prove that every subbrace of a finite direct product of square-free skew braces is complemented. Finally, we extend our investigation to the infinite setting, proving that every multipermutational complemented skew brace is hypercyclic and that all its skew subbraces whose additive groups are $p$-groups are elementary abelian.

\emph{Keywords: Skew left braces, Sylow subbraces, Complements, Multipermutationality}

    \emph{Mathematics Subject Classification (2020):
    16T25, 
    20D10,
	20D20, 
  }
\end{abstract}

\section{Introduction}
The study of complemented subgroups plays an important role in the structural theory of finite groups. Recall that a subgroup $H$ of a group $G$ is said to be \emph{complemented} in $G$ if there exists a subgroup $K$ of $G$ such that $G=HK$ and $H\cap K=1$. Groups in which every subgroup is complemented are called \emph{complemented groups} and form a particularly interesting class, since the existence of complements imposes strong restrictions on the subgroup structure and provides a natural framework for understanding group decompositions and factorisations.

One of the classical results in this direction is due to Hall \cite{Hall37-complemented} and the fourth author and Guo \cite{BallesterGuo99:comp}, establishing the following characterisation of finite complemented groups.

\begin{teo}\label{teo-Hall-complemented}
Let $G$ be a finite group. Then the following statements are equivalent:
\begin{enumerate}
\item $G$ is complemented.
\item Every subgroup of prime order of $G$ is complemented in $G$.
\item $G$ is supersoluble and its Sylow subgroups are elementary abelian.
\item $G$ is isomorphic to a subgroup of a direct product of square-free finite groups.
\end{enumerate}
\end{teo}

The aim of this paper is to investigate the extent to which this characterisation can be generalised to the setting of skew left braces.

A \emph{skew (left) brace} is a set $B$ endowed with two group operations $(B,+)$, the \emph{additive group}, and $(B,\cdot)$, the \emph{multiplicative group}, that satisfy the following distributivity-like property:
\[a(b+c)=ab-a+ac\text{ for all }a,b,c\in B.\]

Skew braces have become an important algebraic tool in the study of set-theoretic solutions of the Yang--Baxter equation. The algebraic structure of a skew brace determines a solution, while properties such as ideals, skew subbraces, quotients and decompositions play an important role in understanding the associated solutions. Consequently, identifying classes of skew braces with particularly strong decomposition properties may provide new approaches to constructing and studying solutions of the Yang--Baxter equation.

One of the key facts about skew braces is that the multiplicative group acts via automorphisms on the additive group by means of the \emph{lambda map}
\[\lambda\colon (B,\cdot)\longrightarrow \Aut(B,+),\quad \lambda_a(b)=-a+ab\text{ for all }a,b\in B.\]
In case this action is trivial, we say that $B$ is a \emph{trivial skew brace}. In this case $a+b=ab$ for all $a,b\in B$. Furthermore, if $B$ is trivial and $(B,+)$ is abelian (respectively elementary abelian) we say that $B$ is \emph{abelian} (respectively \emph{elementary abelian}).

A subset $A$ of $B$ is a \emph{skew subbrace} of $B$ if $A$ is subgroup of $(B,+)$ and $(B,\cdot)$ simultaneously. A skew subbrace of a brace $B$ is called simply \emph{subbrace} throughout.

This leads us to consider skew braces in which every subbrace admits a complement that is itself a subbrace. According to \cite[Definition~8]{BallesterEstebanJimenezPerezC24-solubleskewbraces}, a subbrace $A$ of $B$ is \emph{complemented} in $B$ if there exists a subbrace $C$ such that
\[
B=A+C=AC,\qquad\qquad
A\cap C=0.
\]

In this case, $C$ is called a \emph{complement} of $A$ in $B$. 

We shall refer to skew braces where every subbrace is complemented as \emph{complemented skew braces}.

Recall that a \emph{minimal subbrace} of a brace $B$ is a subbrace of prime order. 

Our first main result establishes a structural consequence of the complementability of minimal subbraces.

\begin{athm}\label{teo-2imp3}
Let $B$ be a finite skew brace in which every minimal subbrace is complemented in $B$. Then $B$ is supersoluble and every subbrace of prime power order is elementary abelian. 
\end{athm}

In particular, this result provides an analogue of the implication from condition 2 to condition 3 in the classical characterisation of complemented groups. Moreover, the complementability of a finite skew brace imposes strong restrictions on both of its underlying group structures.

\begin{cor}\label{cor-2imp3}
Let $B$ be a finite complemented skew brace. Then $(B,+)$ and $(B,\cdot)$ are complemented groups.
\end{cor}

However, the converse of Theorem~\ref{teo-2imp3} does not hold in general. More precisely, we show that the conditions of supersolubility and having elementary abelian $p$-subbraces are not sufficient to guarantee that every minimal subbrace of a finite skew brace is complemented.

\begin{athm}\label{teo-3imp1}
There exists a finite supersoluble skew brace $B$ in which every $p$-subbrace is elementary abelian, but $B$ has minimal subbraces that are not complemented. 
\end{athm}

Thus, unlike the group case, the structural conditions appearing in the classical characterisation do not provide a complete characterisation of complemented skew braces. This illustrates the additional difficulties arising from the compatibility between the additive and multiplicative operations.

We also investigate the role of square-free skew braces, which provide a natural analogue of the groups appearing in condition 4 of the classical result. 

\begin{athm}\label{teo-4imp1}
Every subbrace of a finite direct product of square-free skew braces is complemented.
\end{athm}

Finally, we consider infinite complemented braces. Although the finite characterisation of complemented groups relies on finiteness in an essential way, an infinite analogous structural result can be obtained for multipermutational skew braces.

\begin{athm}\label{teo-inf}
Let $B$ be a multipermutational complemented skew brace. Then $B$ is hypercyclic and every subbrace whose additive group is a $p$-group is elementary abelian.
\end{athm}

\section{Preliminaries}
Let $B$ be skew brace. If $X$ is a subset of $B$ we denote by $\langle X\rangle$ (respectively $\langle X\rangle_+$ and $\langle X\rangle_\bullet$) the smallest subbrace (additive subgroup and multiplicative subgroup respectively) of $B$ that contains $X$.

The \emph{star product} of two elements of $B$ is defined by
\[a\ast b=-a+ab-b=\lambda_a(b)-b\text{ for all }a,b\in B.\]
The star product of two subsets $X,Y\subseteq B$ is defined by
\[X\ast Y=\langle x\ast y\;|\;x\in X,\; y\in Y\rangle_+.\]

A subbrace $A$ of $B$ is a \emph{left ideal} of $B$ if it is $\lambda$-invariant, or equivalently, $B\ast A\subseteq A$. If $A$ is also normal in the additive and multiplicative group we say that $A$ is an \emph{ideal} of $B$. The relevance of ideals lies in the fact that we can make quotients with them. 
Two important subbraces of $B$ that play a key role here are the \emph{kernel of lambda} and the \emph{fixed points set}.
\begin{align*}
\Ker\lambda&=\{a\in B\;|\;\lambda_a(b)-b=a\ast b=0\text{ for all }b\in B\},\\
\Fix(B)&=\{b\in B\;|\;\lambda_a(b)-b=a\ast b=0\text{ for all }a\in B\}.
\end{align*}
The former is a subbrace while the latter is a left ideal of $B$. An important ideal of $B$ is the \emph{socle} of $B$
\begin{align*}
\Soc(B)&=\Ker\lambda\cap Z(B,+).
\end{align*}
We can construct a socle series as follows
\[\Soc_1(B)=\Soc(B),\quad \Soc_{n+1}(B)/\Soc_n(B)=\Soc(B/\Soc_n(B)),\quad n\geq 1.\]
We will say that $B$ is \emph{multipermutational} if there exists $n\geq 1$ such that $\Soc_n(B)=B$. Moreover, the smallest $n\geq 1$ for which this happens is called the \emph{multipermutational level} of $B$.

According to \cite[Definition~18]{BallesterEstebanJimenezPerezC24-solubleskewbraces}, a skew brace $B$ is said to be \emph{soluble} if there exists an ideal series
\[0=I_0\subseteq I_1\subseteq\cdots\subseteq I_{n-1}\subseteq I_n=B\]
of $B$ such that $I_{i+1}/I_i$ is abelian.

We will say that an ideal $I$ of a skew brace $B$ is \emph{$B$-hypercyclic} if there exists an ascending series of ideals of $B$
\[0=I_0\subseteq I_1\subseteq\cdots \subseteq I_\alpha\subseteq I_{\alpha+1}\subseteq \cdots\subseteq I_\gamma=I\]
(here $\alpha$ and $\gamma$ are ordinal numbers) such that $I_{\alpha+1}/I_\alpha$ has prime order. If the series is finite, we will say that $I$ is 
\emph{$B$-supersoluble} (see \cite[Definition~3.1]{BallesterEstebanFerraraPerezCTrombetti24-forumsigma}). If $B=I$, we will simply say that $B$ is \emph{hypercyclic}.

Let $p$ be a prime. A \emph{$p$-subbrace} of $B$ is a subbrace $A$ of $B$ such that $(A,+)$ and $(A,\cdot)$ are $p$-groups, that is, groups in which every element has order a power of $p$. In the finite case that means that $|A|$ is a power of $p$.   

If $B$ is a finite skew brace and $|B| = p^tm$, where $p$ a prime and $(p, m) = 1$, $B$ has subbraces of order $p^t$ (see \cite{ballesterbolinches2026finitetrifactorisedgroupssylow}). Such $p$-subbraces are called \emph{Sylow $p$-subbraces}. 

\section{Proof of Theorem~\ref{teo-2imp3}}

The following result is elementary, therefore, we omit its proof. 

\begin{prop}\label{prop-compl-closed-subbr-quo}
Let $B$ be a complemented skew brace, $A$ a subbrace of $B$ and $I$ an ideal of $B$. Then $A$ and $B/I$ are complemented skew braces.
\end{prop}



\begin{prop}\label{prop-compl-p-skb}
Let $B$ be a finite $p$-skew brace. If every minimal subbrace of $B$ is complemented in $B$, then $B$ is elementary abelian.
\end{prop}
\begin{proof}
Suppose that the result is false, and let $B$ be a minimal counterexample. By Proposition~\ref{prop-compl-closed-subbr-quo} and the minimality of $B$, all subbraces are elementary abelian. We can also suppose that $B$ is a non-trivial brace, because in other case, the result is true by Theorem~\ref{teo-Hall-complemented}.

By \cite[Theorem~4.8 and Corollary~2.27]{CedoSmoktunowiczVendramin19},  $\Fix(B)\neq 0$. Since $\Fix(B)$ is a trivial brace, every subgroup of $\Fix(B)$ is a subbrace of $B$. Hence, we can consider a subbrace $F$ of $B$ of order $p$ contained in $\Fix(B)$. Then there exists a proper subbrace $C$ of $B$ such that $B=F+C = FC$ and $F \cap C = 0$. Let $c_0\in C$ and $b\in B$. Then there exists $f\in F$ and $c\in C$ such that $b=f+c$. Since $f\in\Fix(B)$ and $C$ is elementary abelian, we have that
\[\lambda_{c_0}(b)=\lambda_{c_0}(f+c)=\lambda_{c_0}(f)+\lambda_{c_0}(c)=f+c = b.\]
Hence, $c_0\in\Ker\lambda$, and $C\subseteq\Ker\lambda$. Therefore $|B:\Ker\lambda|\leq p$. Since $B$ is not trivial, we have that $|B:\Ker\lambda|=p$ and $C=\Ker\lambda$. Furthermore, this forces $\Ker\lambda\cap\Fix(B)=0$, because otherwise we could suppose that $F\subseteq\Ker\lambda\cap \Fix(B)$ which is a contradiction. Therefore, $\Fix(B)$ is a complement of $\Ker\lambda$ in $B$ and $|\Fix(B)|=p$.

On the other hand, $\Ker\lambda\neq 0$, since in other case, $B=\Fix(B)$ and $B$ would be trivial. Since $\Ker\lambda$ is a trivial brace, $\Ker\lambda$ contains a subbrace $K$ of order $p$. Let $D$ be a complement of $K$ in $B$. Let $d_0\in D$ and $b\in B$. Then there exists $k\in K$ and $d\in D$ such that $b=kd$. Since $k\in\Ker\lambda$ and $D$ is elementary abelian, we have that
\[\lambda_b(d_0)=\lambda_{kd}(d_0)=\lambda_k(\lambda_d(d_0))=\lambda_k(d_0)=d_0.\]
Therefore, $d_0\in\Fix(B)$ and so $D\subseteq\Fix(B)$.
Analogously to the previous case, $K$ is complemented by $\Fix(B)$ in $B$. Hence $K=\Ker\lambda$ and $|\Ker\lambda|=p$.

Finally, $|B|=|\Ker\lambda||\Fix(B)|=p^2$. In this case $(B,+)$ and $(B,\cdot)$ are abelian groups, therefore, $\Ker\lambda\subseteq \Fix(B)$, a contradiction.\qedhere
\end{proof}

\begin{prop}\label{prop-compl-Syl-tow}
Let $B$ be a finite skew brace with every minimal subbrace complemented in $B$, and let $p$ be the largest prime dividing $|B|$. Then $B$ has a Sylow $p$-subbrace that is an ideal of $B$. 

\end{prop}
\begin{proof}
We argue by induction on $|B|$. If $B$ is a $p$-skew brace, then $B$ is its Sylow $p$-subbrace and the result holds. Suppose that $|B|$ is divisible by at least two primes. By Cauchy's Theorem for braces (see \cite{ballesterbolinches2026finitetrifactorisedgroupssylow}), $B$ has a subbrace $Q$ of order $q$, where $q$ is the smallest prime dividing $|B|$. Let $C$ be a complement of $Q$ in $B$. 

By induction, $C$ has a Sylow $p$-subbrace $P$ that is an ideal of $C$. In particular, $(P,+)$ and $(P,\cdot)$ are characteristic in $(C,+)$ and $(C,\cdot)$ respectively. Since $|B:C|=q$ and $q$ is the smallest prime dividing $|B|$, we have that $(C,+)$ and $(C,\cdot)$ are normal in $(B, +)$ and $(B, \cdot)$ respectively. Hence, $(P,+)$ and $(P,\cdot)$ are normal in $(B,+)$ and $(B,\cdot)$ respectively. Since $q\neq p$, we have that $P$ is a Sylow subbrace of $B$. Since $P$ is characteristic in $(B,+)$ and normal in $(B,\cdot)$, it follows that $P$ is an ideal of $B$.\qedhere
\end{proof}

\begin{proof}[Proof of Theorem~\ref{teo-2imp3}]

We argue by induction on $|B|$. Let $p$ be the largest prime dividing $|B|$, by Proposition~\ref{prop-compl-Syl-tow}, $B$ has a Sylow $p$-subbrace $P$ that is an ideal of $B$. By Proposition~\ref{prop-compl-p-skb}, $P$ is elementary abelian as $P$ is a complemented skew brace by Proposition~\ref{prop-compl-closed-subbr-quo}. Let $T/P$ be a minimal subbrace of $B/P$. Then $|T/P| = q$ for some $q \neq p$. By the Schur-Zassenhaus theorem for skew braces (see \cite{ballesterbolinches2026finitetrifactorisedgroupssylowdominance}), $T$ has a minimal subbrace $S$ such that $T = S + P = SP$ and $S \cap P = 0$. Let $U$ be a complement of $S$ in $B$. Then $B/P = T/P + (U + P)/P = (T/P)((U+P)/P)$ and $T \cap (U + 
P) = (T \cap U) + P = P$ (note that $T\cap U$ is contained in $P$). Therefore, every minimal subbrace of $B/P$ is complemented in $B/P$. By induction $B/P$ is supersoluble, and $P$ is elementary abelian. Therefore, $B$ is soluble. 

Let $N$ be a minimal ideal of $B$. By \cite{BallesterEstebanJimenezPerezC24-solubleskewbraces}, $N$ is abelian. Consider $A$ a subbrace of $N$ of order $p$, and let $C$ be a complement of $A$ in $B$. Since $N$ is an ideal of $B$, $N\cap C$ is an ideal of $C$. Moreover, $N$ is abelian; therefore, $A$ normalises $N\cap C$ additively and multiplicatively. Hence, $N\cap C$ is normal in $C + A =CA = B$ additively and multiplicatively. Let $c_0\in N\cap C$ and $b\in B$. Then there exists $a\in A$ and $c\in C$ such that $b=ac$. Since $N$ is abelian,
\[\lambda_b(c_0)=\lambda_{ca}(c_0)=\lambda_c(\lambda_a(c_0))=\lambda_c(c_0)\in N\cap C.\]
Therefore, $N\cap C$ is an ideal of $B$. The minimality of $N$ yields $N\cap C=0$. It follows that $N=A+(N\cap C)=A$. Hence, $N$ is of order $p$. By induction, $B/N$ is supersoluble, and so $B$ is supersoluble.

Let $q$ be a prime. By Proposition~\ref{prop-compl-closed-subbr-quo}, every minimal subbrace of a $q$-subbrace $S$ of $B$ is complemented. Thus, by Proposition~\ref{prop-compl-p-skb}, $S$ is elementary abelian. \qedhere
\end{proof}

\begin{proof}[Proof of Corollary~\ref{cor-2imp3}]
By Theorem~\ref{teo-2imp3}, $B$ is supersoluble and every subbrace of prime power order is elementary abelian. In particular, the Sylow subbraces of $B$ are elementary abelian. By \cite[Lemma~3.5]{BallesterEstebanFerraraPerezCTrombetti24-forumsigma}, $(B,+)$ and $(B,\cdot)$ are supersoluble. Consequently, both $(B,+)$ and $(B,\cdot)$ are supersoluble groups with elementary abelian Sylow subgroups. By Theorem~\ref{teo-Hall-complemented}, $(B,+)$ and $(B,\cdot)$ are complemented groups.\qedhere
\end{proof}

\section{The example of Theorem~\ref{teo-3imp1}}

According to \cite{BallesterEstebanJimenezPerezC24-solubleskewbraces}, a subbrace $S$ of a skew brace $B$ is called a \emph{maximal subbrace} of $B$ if $S$ is the only proper subbrace of $B$ containing~$S$. The \emph{Frattini} subbrace $\Phi(B)$ of $B$ is defined as the intersection of all maximal subbraces of $B$, if such exist, and $B$ otherwise (\cite[Definition~24]{BallesterEstebanJimenezPerezC24-solubleskewbraces}). If $B$ is finite, every proper subbrace of $B$ is contained in a maximal subbrace of $B$ and so $\Phi(B)$ is a proper subbrace of $B$.

\begin{lemma}\label{lemma-compl-fratt-0}
Let $B$ be a finite skew brace with every minimal subbrace complemented in $B$. Then $\Phi(B)=1$.
\end{lemma}
\begin{proof}
Assume that $\Phi(B) \neq 1$, and let $p$ be a prime dividing $|\Phi(B)|$. Then, by Cauchy's theorem for skew braces (see \cite{ballesterbolinches2026finitetrifactorisedgroupssylow}), $\Phi(B)$ has a subbrace $A$ of order $p$. Let $C$ be a complement of $A$ in $B$. Then  $B = A + C$. Since $C$ is a proper subbrace of $B$, there exists a maximal subbrace $D$ of $B$ such that $C \subseteq D$. Since $A \subseteq D$, it follows that $B = D$, a contradiction. Hence $\Phi(B)=1$.
\end{proof}

\begin{ex}
Let $B$ be the skew brace defined in \cite[Example~3.2]{CampPerez2026Carter}. This skew brace corresponds to \texttt{S(18, 14)} of the \texttt{YangBaxter} library \cite{VendraminKonovalov25-YangBaxter-0.10.7} of \texttt{GAP} \cite{GAP4-15-1}.

Notice that $\Ker\lambda=\langle b, x\rangle_+$ and $\Fix(B)=\langle a,b\rangle_+$.

$\langle b\rangle_+=\langle b\rangle_\bullet$ is normal in $(B,+)$ and $(B,\cdot)$ and lambda invariant. Therefore, it is a cyclic ideal. Moreover, $B/\langle b\rangle_+$ is abelian; thus, $B$ is supersoluble.

Every $2$-subbrace has order $2$, so it is elementary abelian. Also, every $3$-subbrace is elementary abelian since they are contained in $\Fix(B)$, which it is the only Sylow $3$-subgroup of $(B,+)$.

Furthermore, $\Ker\lambda$ and $\Fix(B)$ are maximal subbraces since their indices is prime. They also contain every $2$-subgroup and $3$-subgroup of $(B,+)$ respectively. We show that $\Ker\lambda$ is the only subbrace of order $6$. The subgroups of $(B,+)$ with order $6$ are
\[\langle a,x\rangle_+,\quad \langle a,b+x\rangle_+,\quad \langle a,2b+x\rangle_+,\text{ and } \Ker\lambda=\langle b,x\rangle_+.\]
Since
\[a\cdot x=a+b+x,\quad a\cdot (b+x)=a+2b+x,\text{ and } a\cdot (2b+x)=a+x,\]
the first three subgroups are not closed under the product, hence, $\Ker\lambda$ is the only subbrace of $B$ with order $6$.

In conclusion, $\Ker\lambda$ and $\Fix(B)$ are the only maximal subbraces. Therefore, $\Phi(B)=\Ker\lambda\cap \Fix(B)=\langle b\rangle_+$. By Lemma~\ref{lemma-compl-fratt-0}, $B$ is not complemented.\qedhere
\end{ex}

\section{Proof of Theorem~\ref{teo-4imp1}}

\begin{proof}[Proof of Theorem~\ref{teo-4imp1}]

As in the case of groups, the finite direct product of complemented skew braces is complemented (see \cite{Hall37-complemented}). By Proposition~\ref{prop-compl-closed-subbr-quo}, subbraces of complemented skew braces are complemented. Hence, we only need to show that a square-free skew brace is complemented. Suppose that $B$ is a square-free skew brace. By \cite[Corollary~3.9]{BallesterEstebanFerraraPerezCTrombetti24-forumsigma}, $B$ is supersoluble, in particular it is soluble. Hence, by Hall's theorem for soluble braces (see \cite{ballesterbolinches2026finitetrifactorisedgroupssylow}), $B$ has a subbrace for every integer dividing $|B|$ and, therefore, every subbrace has a complement in $B$.\qedhere
\end{proof}

\section{Proof of Theorem~\ref{teo-inf}}

%
%

\begin{lemma}\label{lemma-finite-gen-subbr}
Let $p$ be a prime and let $B$ be a skew brace such that $(B,+)$ is $p$-elementary abelian. Let $f\in\Fix(B)$ and $s\in\Soc(B)$. Then $\langle f, s\rangle$ is finite.
\end{lemma}
\begin{proof}
Let $n\in\Z$ and write $s_n=\lambda_{nf}(s)=\lambda_{f^n}(s)$. Consider
\[A=\langle f, s_n\;|\;0\leq n\leq p-1\rangle_+=\bigg\{\alpha f+\displaystyle\sum_{n=0}^{p-1}\beta_n s_n\;\bigg|\;0\leq \alpha,\beta_n\leq p-1 \bigg\}.\]
Obviously, $A$ is a subgroup of $(B,+)$. Let us see that is a subbrace. First, notice that $s_n=s_m$ if $n\equiv m\mod p$, hence, $s_n\in A$ for all $n\in\Z$. Moreover, $s_n=nf\ast s+s\in\Soc(B)$ for all $n\in\Z$. Consider $a,b\in A$. hen 
\[a=\alpha f+\sum_{n=0}^{p-1}\beta_n s_n=\bigg(\sum_{n=0}^{p-1}\beta_n s_n\bigg)(\alpha f),\quad b=\alpha' f+\sum_{n=0}^{p-1}\beta_n' s_n\]
for some $0\leq\alpha,\alpha',\beta_n,\beta_n'\leq p-1$. Then
\begin{align*}
ab&=a+\lambda_a(b)=\lambda_{\alpha f}\bigg(\alpha' f+\sum_{n=0}^{p-1}\beta_n' s_n \bigg)=\alpha' f+\sum_{n=0}^{p-1}\beta_n'\lambda_{\alpha f}(s_n)
\end{align*}
Notice that
\[\lambda_{\alpha f}(s_n)=\lambda_{f^\alpha}(\lambda_{f^n}(s))=\lambda_{f^{(\alpha + n)}}(s)=\lambda_{(\alpha+ n)f}(s)=s_{\alpha+n}\in A.\]
Hence, $ab\in A$. Moreover,
\begin{align*}
a^{-1}&=(\alpha f)^{-1}\bigg(\sum_{n=0}^{p-1}\beta_n s_n\bigg)^{-1}=(-\alpha f)\bigg(-\sum_{n=0}^{p-1}\beta_n s_n\bigg)\\&=((-\alpha) f)\bigg(\sum_{n=0}^{p-1}(-\beta_n) s_n\bigg).
\end{align*}
Since $f$ and $s_n$ have additive orders either $1$ or $p$, we have that $(-\alpha)f$ and $\sum_{n=0}^{p-1}(-\beta_n) s_n$ are in $A$. Since the product of any two elements of $A$ belongs to $A$, it follows that $a^{-1}\in A$. Hence $A$ is a finite subbrace of $B$ containing $\langle f,s\rangle$.\qedhere
\end{proof}

\begin{prop}\label{prop-inf-compl-multiperm-p-brace-elem-ab}
Let $B$ be a multipermutational complemented skew brace such that $(B,+)$ is a $p$-group for a prime $p$. Then $B$ is elementary abelian.
\end{prop}
\begin{proof}
First, let us prove that $(B,+)$ is elementary abelian. By \cite[Theorem~2.20]{CedoSmoktunowiczVendramin19}, $(B,+)$ is nilpotent. We argue by induction on the nilpotency class of $(B,+)$.

Suppose that $(B,+)$ is abelian, and consider $A=\{a\in B\;|\; pa=0\}$. Suppose that $A\neq B$. Since $(B,+)$ is abelian, $(A,+)$ is a characteristic subgroup of $(B,+)$; therefore, $A$ is a subbrace of $B$. Let $C$ be a subbrace of $B$ such that $B=A+C$ and $A\cap C=0$. If $C \neq 0$, then $C$ has an element of order $p$ by Cauchy's theorem, a contradiction. Hence $C = 0$, $B=A$ and $(B,+)$ is elementary abelian.

Now suppose that the nilpotency class of $(B,+)$ is $n$ and that the result is true for multipermutational complemented skew braces whose additive group has nilpotency class $n-1$. Since $Z(B,+)$ is a characteristic subgroup of $(B,+)$, it is a subbrace of $B$, so that there exists a subbrace $D$ complementing $Z(B,+)$ in $B$. We have that $(D,+)\cong (B,+)/Z(B,+)$. By \cite[Lemma~2.6 and Theorem~2.20]{CedoSmoktunowiczVendramin19} and Proposition~\ref{prop-compl-closed-subbr-quo}, $D$ is a multipermutational complemented skew brace and the nilpotency class of $(D,+)$ is less than $n$. By induction, $(D,+)$ is elementary abelian. Since $Z(B,+)$ is an abelian subbrace of $B$, we have that $Z(B,+)$ is elementary abelian. Finally, $(B,+)=Z(B,+)+(D,+)$, so that $(B,+)$ is elementary abelian.

Next, we show that $B$ is a trivial brace by induction on the multipermutational level $n$ of $B$. It is clear that we may assume $n > 1$ and that the result is true for multipermutational complemented skew braces of multipermutational level $n-1$. Let $C$ be a complement of $\Soc(B)$ in $B$. Then $C\cong B/\Soc(B)$. By \cite[Lemma~2.6 and Theorem~2.20]{CedoSmoktunowiczVendramin19} and Proposition~\ref{prop-compl-closed-subbr-quo}, $C$ is a multipermutational complemented skew brace with multipermutational level $n-1$. By induction, $C$ is a trivial brace. Let $c_0\in C$ and $b\in B$. Then there exists $s\in \Soc(B)$ and $c\in C$ such that $b=sc$. Since $s\in\Soc(B)$ and $C$ is trivial, we have that
\[\lambda_{b}(c_0)=\lambda_{sc}(c_0)=\lambda_{s}(\lambda_c(c_0))=c_0.\]
Hence, $c_0\in\Fix(B)$, and $C\subseteq\Fix(B)$.

Consider $f\in\Fix(B)$ and $s\in\Soc(B)$. By Lemma~\ref{lemma-finite-gen-subbr}, $A=\langle f,s\rangle$ is finite. Therefore, $A$ is a finite complemented $p$-skew brace. Thus, by Theorem~\ref{teo-2imp3}, $A$ is trivial, in particular, $\lambda_f(s)=s$.

Let $a,b\in B$. Since $B=\Fix(B)+\Soc(B)=\Fix(B)\Soc(B)$, there exists $f,f'\in\Fix(B)$ and $s,s'\in\Soc(B)$ such that $a=fs$ and $b=f'+s'$. Therefore,
\[\lambda_a(b)=\lambda_{fs}(f'+s')=f'+\lambda_f(s')=f'+s'=b.\]
Consequently, $B$ is a trivial brace.\qedhere
\end{proof}

Recall that if $G$ is a periodic abelian group and $p$ is a prime, then $P = \{g \in G : p^ng = 0$ for some $n \geq 0$\} is the unique Sylow $p$-subgroup of $G$. If $G$ is a skew brace, then $P$ is the Sylow $p$-subbrace of $G$.

\begin{lemma}\label{lemma-Syl-id-soc}
Let $B$ be a skew brace and let $I$ be an ideal of $B$ contained in $\Soc(B)$. If $\Soc(B)$ is periodic, then the Sylow subgroups of $(I,+)$ are ideals of $B$.
\end{lemma}
\begin{proof}
 Let $(P,+)$ be a Sylow $p$-subbrace of $(I,+)$ for some prime $p$. Since $I$ is a trivial brace, $P$ is a subbrace of $B$. Moreover, $(P,+)$ and $(P,\cdot)$ are characteristic in $(I,+)$ and $(I,\cdot)$, respectively. Since $(I,+)$ and $(I,\cdot)$ are normal in $(B,+)$ and $(B,\cdot)$, respectively, it follows that $(P,+)$ and $(P,\cdot)$ are characteristic in $(B,+)$ and $(B,\cdot)$ respectively. Finally, since $I$ is $\lambda$-invariant and $(P,+)$ is characteristic in $(I,+)$,  we have that $P$ is $\lambda$-invariant. Therefore, $P$ is an ideal of $B$.\qedhere
\end{proof}

\begin{lemma}\label{l}
If $B$ is a complemented skew brace, then $\Soc(B)$ is periodic.
\end{lemma}

\begin{proof}
Assume that $\Soc(B)$ has an element $s$ of infinite additive 
order. Since $\Soc(B)$ is an abelian subbrace, the 
additive infinite cyclic group $\langle s \rangle$ and its subgroup $\langle 2s \rangle$ are subbraces of $B$. By Proposition~\ref{prop-compl-closed-subbr-quo}, $\langle 2s \rangle$ has a complement in $\langle s \rangle$, a contradiction. Hence, $\Soc(B)$ is periodic.\qedhere
\end{proof}

\begin{prop}\label{prop-ex-id-prime-ord}
Let $B$ be a complemented skew brace and let $I$ be a non-zero ideal of $B$ contained in $\Soc(B)$. Then there exists an ideal of $B$ of prime order contained in $I$.
\end{prop}
\begin{proof}
By Lemma~\ref{l}, $I$ is a periodic abelian skew brace. Let $(P,+)$ be a non-zero Sylow subgroup of $(I,+)$. By Lemma~\ref{lemma-Syl-id-soc}, $P$ is an ideal of $B$. By Proposition~\ref{prop-compl-closed-subbr-quo}, $P$ is complemented. Since $P$ is a trivial brace, we can consider a subbrace of $P$ of order $p$, and then take a complement $U$ in $P$. Then $U$ is a subbrace of $P$ such that $|P:U|=p$.

Consider a complement $C$ of $U$ in $B$, and let $V=C\cap P$. Then $V\cap U\subseteq U\cap C=0$, and
\[P=P\cap B=P\cap (UC)=U(P\cap C)=UV.\]
Hence, $|V|=p$. By Proposition~\ref{prop-compl-closed-subbr-quo}, $C$ is complemented. Let $Q$ be a complement of $V$ in $C$. Since $P$ is an ideal of $B$, we have that $V=C\cap P$ is an ideal of $C$; in particular, $V$ is lambda invariant by the elements of $Q$. Moreover, since $V\subseteq P$, we have that $V$ is lambda invariant by the elements of $P$. Thus, $V$ is lambda invariant by the elements of $PQ=UVQ=UC=B$. Consequently, $V$ is a left ideal of $B$. Finally, note that every left ideal of $B$ contained in $\Soc(B)$ is an ideal.\qedhere
\end{proof}

\begin{lemma}\label{lemma-B-hypercyc-quo}
Let $B$ be a skew brace and let $I\subseteq J$ be ideals of $B$ such that $J$ is $B$-hypercyclic. Then $J/I$ is $B/I$-hypercyclic.
\end{lemma}
\begin{proof}
Let
\[J_0=0\subseteq\cdots\subseteq J_\alpha\subseteq J_{\alpha+1}\subseteq\cdots\subseteq J_\gamma=J\]
be an ascending series of ideals of $B$ such that $J_{\alpha+1}/J_\alpha$ has prime order. Denote $I_\alpha=(J_\alpha+I/I)$. Then
\[I_0=0\subseteq\cdots\subseteq I_\alpha\subseteq I_{\alpha+1}\subseteq\cdots\subseteq I_\gamma=J/I\]
is an ascending series of ideals of $B/I$. Moreover,
\begin{align*}
I_{\alpha+1}/I_\alpha &\cong (J_{\alpha+1}+I)/(J_\alpha+I)=(J_{\alpha+1}+(J_\alpha+I))/(J_\alpha+I)\\&\cong J_{\alpha+1}/(J_{\alpha+1}\cap (J_\alpha+I))=J_{\alpha+1}/(J_\alpha+(J_{\alpha+1}\cap I)).
\end{align*}
Since $J_{\alpha+1}/J_\alpha$ has prime order, either $I_{\alpha+1}/I_\alpha$ has prime order or $I_{\alpha+1}=I_\alpha$. Therefore, if we remove the repeated terms from the chain we have an ascending chain of ideals of $B/I$ whose successive quotients have prime order. Thus, $J/I$ is $B/I$-hypercyclic.\qedhere
\end{proof}

\begin{prop}\label{prop-soc-B-hypercyc}
Let $B$ be a complemented skew brace. Then $\Soc(B)$ is $B$-hypercyclic.
\end{prop}
\begin{proof}
Let $\mathfrak{H}$ be the set of all $B$-hypercyclic ideals of $B$ contained in $\Soc(B)$. Obviously $0\in \mathfrak{H}$, so it is not an empty set. Let
\[I_0=0\subseteq\cdots\subseteq I_\alpha\subseteq I_{\alpha+1}\subseteq\cdots\]
be an ascending chain of elements of $\mathfrak{H}$. Consider $I=\cup_{\alpha}I_\alpha$. Clearly, $I$ is an ideal of $B$ contained in $\Soc(B)$. Moreover, by Lemma~\ref{lemma-B-hypercyc-quo}, we can find an ascending series of ideals of $B$ from $I_\alpha$ to $I_{\alpha+1}$ such that the quotients have prime order. Therefore, we can construct an ascending series of ideals of $B$ from $0$ to $I$ such that the quotients have prime order. Therefore, $I\in\mathfrak{H}$.

By Zorn's lemma, there exists $M$, a maximal element of $\mathfrak{H}$. If $M\subsetneq \Soc(B)$, then $\Soc(B)/M$ is a non-zero ideal of $B/M$ contained in $\Soc(B/M)$. By Proposition~\ref{prop-compl-closed-subbr-quo}, $B/M$ is a complemented skew brace. Then, by Proposition~\ref{prop-ex-id-prime-ord}, there exists an ideal $I/M$ of $B/M$ of prime order contained in $\Soc(B)/M$. Since $M$ is $B$-hypercyclic and $I/M$ has prime order, $I$ is $B$-hypercyclic. This contradicts the maximality of $M$. Thus, $\Soc(B)=M$ is $B$-hypercyclic.\qedhere
\end{proof}

\begin{proof}[Proof of Theorem~\ref{teo-inf}]

We show that $B$ is hypercyclic by induction on the multipermutational level of $B$. By Proposition~\ref{prop-soc-B-hypercyc}, $\Soc(B)$ is $B$-hypercyclic and, by Proposition~\ref{prop-compl-closed-subbr-quo}, $B/\Soc(B)$ is a multipermutational complemented skew brace of smaller multipermutational level. By induction, $B/\Soc(B)$ is hypercyclic, and hence $B$ is hypercyclic.

By \cite[Lemma~2.6 and Theorem~2.20]{CedoSmoktunowiczVendramin19} and Proposition~\ref{prop-compl-closed-subbr-quo}, every subbrace of $B$ is a complemented multipermutational skew brace.  Hence, by Proposition~\ref{prop-inf-compl-multiperm-p-brace-elem-ab}, every subbrace whose additive group is a $p$-group is elementary abelian for every prime $p$.\qedhere
\end{proof}

\subsection*{Acknowledgements}
\begin{sloppypar}
The first author was supported by Guangdong Basic and Applied Research Foundation (2025A1515012072). The second and fourth authors are supported by the grant PID2024-159495NB-I00, funded byMICIU/AEI/10.13039/501100011033 and by ERDF/EU, and by the grant CIAICO/2023/007 from the Conselleria d’Educació, Cultura, Universitats i Ocupació of the Generalitat Valenciana. The second author is funded by the Spanish Ministry of Science, Innovation and Universities through an FPU predoctoral grant FPU24/01319. The third author is grateful to the Universitat de València for their support to host researchers affected by the war in Ukraine.
\end{sloppypar}

\bibliographystyle{plain}
\bibliography{bibgroup}

@article{BallesterEstebanFerraraPerezCTrombetti24-forumsigma,
  title   = {Finite skew braces of square-free order and supersolubility},
  author  = {A. Ballester-Bolinches and R. Esteban-Romero and
             M. Ferrara and V. P{\'e}rez-Calabuig and M. Trombetti},
  journal = {Forum Math. Sigma},
  volume  = {12},
  pages   = {e39},
  year    = {2024},
  doi     = {10.1017/fms.2024.29}
}

@article{BallesterEstebanJimenezPerezC24-solubleskewbraces,
  title   = {Soluble skew left braces and soluble solutions of the {Y}ang--{Baxter} equation},
  author  = {A. Ballester-Bolinches and R. Esteban-Romero and
             P. Jim{\'e}nez-Seral and V. P{\'e}rez-Calabuig},
  journal = {Adv. Math.},
  volume  = {455},
  pages   = {109880},
  year    = {2024},
  doi     = {10.1016/j.aim.2024.109880}
}

@article{BallesterGuo99:comp,
  title   = {On complemented subgroups of finite groups},
  author  = {A. Ballester-Bolinches and X. Guo},
  journal = {Arch. Math. (Basel)},
  volume  = {72},
  number  = {3},
  pages   = {161--166},
  year    = {1999},
  doi     = {10.1007/s000130050317}
}

@misc{CampPerez2026Carter,
  title   = {Carter-like theorems for skew braces},
  author  = {Sergio Camp-Mora and P. P{\'e}rez-Altarriba},
  year    = {2026},
  doi     = {10.48550/arXiv.2609.39430},
  note    = {preprint}
}

@article{CedoSmoktunowiczVendramin19,
  title   = {Skew left braces of nilpotent type},
  author  = {F. Ced{\'o} and A. Smoktunowicz and L. Vendramin},
  journal = {Proc. London Math. Soc.},
  volume  = {118},
  number  = {6},
  pages   = {1367--1392},
  year    = {2019},
  doi     = {10.1112/plms.12209}
}

@manual{GAP4-15-1,
  title        = {{GAP -- Groups, Algorithms, and Programming,
  				    Version 4.15.1}},
  author       = {{The GAP Group}},
  organization = {The GAP Group},
  year         = {2025},
  url          = {https://www.gap-system.org},
  note         = {{\url{http://www.gap-system.org}}}
}

@article{Hall37-complemented,
  title   = {Complemented {G}roups},
  author  = {P. Hall},
  journal = {J. London Math. Soc.},
  volume  = {s1-12},
  number  = {3},
  pages   = {201--204},
  year    = {1937},
  doi     = {10.1112/jlms/s1-12.2.201}
}

@article{ballesterbolinches2026finitetrifactorisedgroupssylow,
  title   = {On finite trifactorised groups and {S}ylow and {H}all theorems for skew braces},
  author  = {A. Ballester-Bolinches and P. P{\'e}rez-Altarriba and
             V. P{\'e}rez-Calabuig},
  journal = {Ricerche Mat.},
  year    = {2026},
  doi     = {10.1007/s11587-026-01167-9}
}

@manual{VendraminKonovalov25-YangBaxter-0.10.7,
  title  = {Yang{B}axter: combinatorial solutions for the
            {Y}ang--{Baxter} equation},
  author = {L. Vendramin and O. Konovalov},
  year   = {2025},
  note   =	{version 0.10.7, \url{https://gap-packages.github.io/YangBaxter/}},
  url    = {https://gap-packages.github.io/YangBaxter/}
}

@misc{ballesterbolinches2026finitetrifactorisedgroupssylowdominance,
  title    = {On finite trifactorised groups and dominance {S}ylow and {H}all theorems for finite skew braces},
  author   = {A. Ballester-Bolinches and P. P{\'e}rez-Altarriba and V. P{\'e}rez-Calabuig},
  year     = {2026},
  note     = {Preprint}
}

\end{document}